\documentclass[graybox]{svmult}

\usepackage{type1cm}        
\usepackage{makeidx}         
\usepackage{graphicx}        
\usepackage{multicol}        
\usepackage[bottom]{footmisc}
\usepackage[numbers]{natbib}

\usepackage{newtxtext}       %
\usepackage[varvw]{newtxmath}       

\usepackage{amsmath,graphicx,textcomp,hyperref,tabularx,xcolor,theorem}

\hypersetup{
	colorlinks=true, 
	linkcolor=blue, 
	urlcolor=red, 
	citecolor=magenta,
	linktoc=all 
}

\usepackage{algorithmicx}
\usepackage{algorithm}
\usepackage{algpseudocode}

\usepackage{xpatch}
\makeatletter
\xpatchcmd{\algorithmic}
{\ALG@tlm\z@}{\leftmargin\z@\ALG@tlm\z@}
{}{}
\makeatother

\newcommand{\backassign}{=:}

\newcommand{\nocomma}{}

\newcommand{\tmmathbf}[1]{\ensuremath{\boldsymbol{#1}}}
\newcommand{\tmop}[1]{\ensuremath{\operatorname{#1}}}

\newcommand{\tmtextit}[1]{\text{{\itshape{#1}}}}
\newcommand{\tmtextrm}[1]{\text{{\rmfamily{#1}}}}

\newcommand{\pdv}[2]{\frac{\partial #1}{\partial #2}}

\newcommand{\ba}{\tmmathbf{a}}
\newcommand{\bb}{\tmmathbf{b}}

\newcommand{\be}{\tmmathbf{e}}
\newcommand{\bc}{\tmmathbf{c}}

\newcommand{\bn}{\tmmathbf{n}}

\newcommand{\bp}{\tmmathbf{p}}

\newcommand{\bq}{\ensuremath{\tmmathbf{q}}}
\newcommand{\bzero}{\tmmathbf{0}}

\newcommand{\bv}{\tmmathbf{v}}
\newcommand{\bx}{\tmmathbf{x}}

\newcommand{\ud}{\text{\tmtextrm{d}}}

\newcommand{\uu}{\tmmathbf{u}}
\newcommand{\uU}{\tmmathbf{U}}

\newcommand{\pf}{\tmmathbf{f}}

\newcommand{\F}{\ensuremath{\tmmathbf{F}}}

\newcommand{\tf}{\tilde{\pf}}
\newcommand{\tF}{\tilde{\F}}

\newcommand{\Div}{\nabla_{\tmmathbf{x}} \cdot}

\newcommand{\xii}{\xi^i}

\newcommand{\polyP}{P}

\newcommand{\ijkc}{i, j, k}

\newcommand{\Oo}{\Omega_o}
\newcommand{\Oip}{\partial \Omega_{o, i}^R}
\newcommand{\Ois}{\partial \Omega_{o, i}^s}
\newcommand{\Oim}{\partial \Omega_{o, i}^L}

\newcommand{\bnr}{\widehat{\tmmathbf{n}}}

\newcommand{\Nnd}{\mathbb{N}_N^d}

\newcommand{\uebp}{\tmmathbf{u}_{e, \tmmathbf{p}}}

\newcommand{\paragraphtoc}[1]{}
\providecommand{\subindex}{=}

\newcommand{\fa}{\tmmathbf{f}^a}
\newcommand{\fv}{\tmmathbf{f}^v}
\newcommand{\tfa}{\tilde{\tmmathbf{f}}^a}
\newcommand{\tfv}{\tilde{\tmmathbf{f}}^v}
\newcommand{\tfs}{\tilde{\tmmathbf{f}}^\alpha}
\newcommand{\tfad}{\tilde{\tmmathbf{f}}^{a \delta}}
\newcommand{\tfvd}{\tilde{\tmmathbf{f}}^{v \delta}}

\newcommand{\tFa}{\tilde{\tmmathbf{F}}^a}
\newcommand{\tFv}{\tilde{\tmmathbf{F}}^v}
\newcommand{\tFs}{\tilde{\tmmathbf{F}}^\alpha}
\newcommand{\tFad}{\tilde{\tmmathbf{F}}^{a \delta}}
\newcommand{\tFvd}{\tilde{\tmmathbf{F}}^{v \delta}}
\newcommand{\vxi}{\ensuremath{\tmmathbf{\xi}}}
\newcommand{\M}{\mathcal{M}}
\newcommand{\mB}{\mathcal{B}}
\newcommand{\Pra}{\text{\tmtextit{Pr}}}
\newcommand{\Grad}{\nabla_{\tmmathbf{x}}}
\usepackage{url}

\makeindex             

\begin{document}
	
\title*{Lax-Wendroff Flux Reconstruction for advection-diffusion equations with error-based time stepping}
\titlerunning{LWFR for advection diffusion equations}
\author{Arpit Babbar\orcidID{0000-0002-9453-370X} and\\ Praveen Chandrashekar\orcidID{0000-0003-1903-4107}}
\authorrunning{Babbar, Chandrashekar}
\institute{Arpit Babbar \at TIFR-CAM, Bangalore - 560065, \email{arpit@tifrbng.res.in}
\and Praveen Chandrashekar \at TIFR-CAM, Bangalore - 560065, \email{praveen@math.tifrbng.res.in}}
\maketitle
\abstract{This work introduces an extension of the high order, single stage Lax-Wendroff Flux Reconstruction (LWFR) of Babbar et al., JCP (2022) to solve second order time-dependent partial differential equations in conservative form on curvilinear meshes. The method uses BR1 scheme to reduce the system to first order so that the earlier LWFR scheme can be applied. The work makes use of the embedded error-based time stepping introduced in Babbar, Chandrashekar (2024) which becomes particularly relevant in the absence of CFL stability limit for parabolic equations. The scheme is verified to show optimal order convergence and validated with transonic flow over airfoil and unsteady flow over cylinder.}
\section{Introduction}
Spectral element methods like Flux Reconstruction (FR)~\cite{Huynh2007} and Discontinuous Galerkin (DG)~\cite{cockburn2000} are very successful in resolving advective flows~\cite{witherden2014,Ranocha2022}. Their first extension to second order PDE like Navier-Stokes is the BR1 scheme~\cite{Bassi1997}. The BR1 scheme retains the superior properties of FR / DG and is even applicable in underresolved turbulent simulations~\cite{Gassner2012}; BR1 was proven to be stable in~\cite{Gassner_Winters_Hindenlang_Kopriva_2018}.

Lax-Wendroff Flux Reconstruction (LWFR) is a spectral element solver for time dependent PDE that performs high order evolution to the next time level in a single stage, minimizing interelement communication which is often the bottleneck for the modern memory bandwidth limited hardware. The general idea of high order LW was introduced in~\cite{Qiu2003} and is based on performing a Taylor's expansion in time and Cauchy-Kowaleski procedure on the solution. This work extends the LWFR scheme of~\cite{babbar2022,babbar2024curved} to second order parabolic equations on curvilinear meshes by making use of the BR1 scheme. LWFR has been extended to curved meshes in~\cite{babbar2024curved} along with the first embedded error-based time stepping method for single stage methods, which was relevant there since a Fourier $L^2$ stability analysis only applies to Cartesian meshes. The controller used for error-based time stepping in~\cite{babbar2024curved} was taken from~\cite{Ranocha2021} where optimal controllers were found for error-based time stepping methods for multistage Runge-Kutta methods, and the scheme was applied to various CFD problems. This work uses the error based time stepping of~\cite{babbar2024curved}, which is especially relevant here since a Fourier CFL stability limit of LWFR is also not known for second order PDE.

This is the first extension of a high order Lax-Wendroff to second order equations. The ADER schemes, which are another class of single stage solvers have also been extended to solve second order PDE in~\cite{Fambri2017} by including additional diffusion in the numerical flux in contrast to the BR1 scheme used here. This is the first work where any single stage method has been combined with the BR1 scheme.

The rest of this paper is organized as follows. In Section~\ref{sec:transformations.parabolic}, the notations and transformations of second order equations from curved element to a reference cube are reviewed. In Section~\ref{sec:lwfr.2nd.order}, the LWFR scheme for second order equations is introduced. The numerical results are shown in Section~\ref{sec:results} and conclusions from the work made in Section~\ref{sec:conclusion}.
\section{Curvilinear coordinates for parabolic equations} \label{sec:transformations.parabolic}
We explain our numerical scheme for a general system of parabolic equations in conservative form in dimension $d$
\begin{equation}
\partial_t  \uu + \Div \fa \left( \uu \right) = \Div \fv \left( \uu, \Grad 
\uu \right) \label{eq:parabolic.eqn}
\end{equation}
with some initial and boundary conditions. Here, $\uu \in \mathbb{R}^p$ is the
solution vector and its gradient is the matrix $\Grad  \uu = \left( \partial_{x_1}  \uu, \ldots,
\partial_{x_d}  \uu \right) \in \mathbb{R}^{p \times d}$. The $\fa, \fv$ are
the advective and viscous fluxes and can be seen as matrices $\pf = \left(
\pf_1, \ldots, \pf_d \right) \in \mathbb{R}^{p \times d}$, $\bx$ is in domain
$\Omega \subset \mathbb{R}^d$ and divergence of a flux is given by $\Div \pf = \sum_{i = 1}^d \partial_{x_i}  \pf_i$. In general, following the notation of~\cite{Gassner_Winters_Hindenlang_Kopriva_2018}, the
action of a vector $\bb \in \mathbb{R}^d$ on $\bv \in \mathbb{R}^p$ gives $\bb
\bv \in \mathbb{R}^{p \times d}$ which is defined component-vice as
\begin{equation}
\bb  \bv = \left( b_j  \bv \right)_{j = 1}^d \label{eq:vector.vector}
\end{equation}
Further, the action of a matrix $\mB = \left( \bb_1, \ldots, \bb_d \right) \in
\mathbb{R}^{d \times d}$ on $\bv = (\bv_1, \ldots, \bv_d)\in \mathbb{R}^{p \times d}$ also gives $\mB 
\bv \in \mathbb{R}^{p \times d}$ defined as
\begin{equation}
	\mathcal{B}  \bv = \sum_{i = 1}^d \bb_i  \bv_i, \qquad \bb_i  \bv_i = (b_{i j}  \bv_i)_{j = 1}^d \label{eq:matrix.operator.notation}
\end{equation}
The second order system~\eqref{eq:parabolic.eqn} will be reduced to a first order system of
equations following the BR1 scheme~{\cite{Bassi1997}} where an auxiliary
variable $\bq$ is introduced
\begin{equation}
\begin{split}
\bq - \Grad  \uu & = \bzero \\
\partial_t  \uu + \Div \pf \left( \uu, \bq \right) & = \bzero
\end{split}\label{eq:parabolic.in.first.order}
\end{equation}
We partition the domain $\Omega$ into non-overlapping quadrilateral elements $\{\Omega_e\}$. The reference element is the cube $\Oo = [- 1, 1]^d$ with coordinates $\vxi = ( \xii )_{i = 1}^d$
and the Cartesian orthonormal basis vectors denoted as $\left\{ \be_i
\right\}_{i = 1}^d$. Corresponding to each element $\Omega_e$, there is a
bijective reference map $\Theta_e$ from the reference cube $\Oo$ onto the
element $\Omega_e$ so that the physical variable is $\bx = \Theta_e( \vxi )$; the element index $e$ will usually be suppressed. We will denote a $d$-dimensional
multi-index as $\bp = (p_i)_{i = 1}^d$. In this work, the reference mapping is
defined using Lagrange interpolation
\begin{equation}
\Theta \left( \vxi \right) = \sum_{\bp \in \Nnd} \widehat{\bx}_{\bp} \ell_{\bp} \left( \vxi \right) \label{eq:ref.map}
\end{equation}
where $\Nnd = \left\{ \bp : p_i \in \{ 0, 1, \ldots, N \}, 1 \leq i \leq d \right\}$ and $\{ \ell_{\bp} \}_{\bp \in \Nnd}$ is the degree $N$ Lagrange polynomial basis corresponding to the points $\{ \vxi_{\bp} \}_{\bp \in \Nnd}$ so that $\Theta ( \vxi_{\bp} ) = \widehat{\bx}_{\bp}$ for all $\bp \in \Nnd$. Thus, the points $\{ \vxi_{\bp} \}_{\bp \in \Nnd}$ are where the reference map will be specified and they will be taken to be the solution points of the Flux Reconstruction scheme throughout this work. In this work, we only use Gauss-Legendre-Lobatto (GLL) solution points which include element interfaces, which is why~\eqref{eq:ref.map} gives a globally continuous reference map.
The basis $\{ \ell_{\bp} \}_{\bp \in \Nnd}$ can be written in terms of the 1-D Lagrange basis $\{ \ell_{p_i} \}_{p_i \subindex 0}^N$ of degree $N$ corresponding to points $\{ \xi_{p_i} \}_{p_i = 0}^N$ using tensor product
polynomials
\begin{equation}
\ell_{\bp} \left( \vxi \right) = \prod_{i = 1}^d \ell_{p_i} \left( \xii \right)\label{eq:lagrange.polynomials}
\end{equation}
The covariant and contravariant basis vectors are used in specifying transformation of~\eqref{eq:parabolic.eqn} to equations in the reference element 
$\Oo$ and are defined as follows.
\begin{definition}[Covariant basis]\label{defn:covariant.basis}
The coordinate basis vectors $\left\{ \ba_i \right\}_{i = 1}^d$ are defined
so that $\ba_i, \ba_j$ are tangent to $\{ \xi^k = \tmop{const} \}$ where $i,
j, k$ are cyclic. They are explicitly given as
\begin{equation}
	\ba_i = \pdv{\bx}{\xii}, \qquad 1 \leq i \leq d
\end{equation}
\end{definition}
\begin{definition}[Contravariant basis]\label{defn:contravariant.basis}
The contravariant basis vectors $\left\{ \ba^i \right\}_{i = 1}^d$ are the
respective normal vectors to the coordinate planes $\{ \xi_i = \tmop{const}
\}_{i = 1}^3$. They are explicitly given as
\begin{equation}
\ba^i = \nabla_{\bx}  \xii, \qquad 1 \leq i \leq d
\end{equation}
\end{definition}
The following relation is satisfied by the contravariant and covariant
vectors~{\cite{Kopriva2006,kopriva2009}}
\begin{equation}
	J \ba^i = J \nabla \xii = \ba_j \times \ba_k
	\label{eq:contravariant.identity}
\end{equation}
where $\left( \ijkc \right)$ are cyclic, and $J$ denotes the Jacobian of the
transformation which also satisfies
\[ J = \ba_i \cdot \left( \ba_j \times \ba_k \right) \qquad (i, j, k) \text{
	cyclic} \]
The divergence of a flux under transformed coordinates can be computed
as~{\cite{Kopriva2006,kopriva2009}}
\begin{equation}
	\Div \pf = \frac{1}{J}  \sum_{i = 1}^d \pdv{}{\xii}  \left( J \ba^i \cdot
	\pf \right) \label{eq:divergence.transform.contravariant}
\end{equation}
Consequently, the gradient of a scalar $\phi$ becomes
\begin{equation}
	\Grad \phi = \frac{1}{J}  \sum_{i = 1}^d \pdv{}{\xii}  \left[ \left( J \ba^i
	\right) \phi \right] \label{eq:scalar.gradient.transform.contravariant}
\end{equation}
The $\left\{ J \ba^i \right\}_{i = 1}^3$ are called the metric terms. Since
the divergence / gradient of a constant function is
zero,~\eqref{eq:contravariant.identity} gives us the metric identity
\begin{equation}
	\sum_{i = 1}^d \pdv{J \ba^i}{\xii} = \bzero
	\label{eq:metric.identity.contravariant}
\end{equation}
Using Leibniz product rule and the metric identity, we get the non-conservative form of gradient of a
scalar $\phi$ and thus of a vector $\uu$ in vector action
notation~\eqref{eq:vector.vector}
\begin{equation}
\Grad \phi = \frac{1}{J}  \sum_{i = 1}^d J \ba^i  \pdv{\phi}{\xii}, \qquad
\Grad  \uu = \frac{1}{J}  \sum_{i = 1}^d J \ba^i  \pdv{\uu}{\xii}
\end{equation}
Following the notation of~{\cite{Gassner_Winters_Hindenlang_Kopriva_2018}},
define the transformation matrix $\M = \left( J \ba^1, \ldots, J \ba^d \right)
\in \mathbb{R}^{d \times d}$ so with matrix action
notation~\eqref{eq:matrix.operator.notation}
\begin{equation}
\Grad  \uu = \frac{1}{J}  \M \nabla_{\vxi}  \uu  \label{eq:grad.transform}
\end{equation}
Within each element $\Omega_e$, performing change of variables with the
reference map $\Theta_e$~(\ref{eq:grad.transform},\ref{eq:divergence.transform.contravariant}), the
first order system~\eqref{eq:parabolic.in.first.order} transforms into
\begin{equation}
\begin{split}
J \bq - \M \nabla_{\vxi}  \uu & = \bzero\\
J \partial_t  \uu + \nabla_{\vxi} \cdot \tfa & = \nabla_{\vxi} \cdot \tfv
\left( \uu, \bq \right)
\end{split} \label{eq:transformed.fo}
\end{equation}
\section{Lax-Wendroff flux reconstruction} \label{sec:lwfr.2nd.order}
The solution of the conservation law will be approximated by piecewise
polynomial functions which are allowed to be discontinuous across the
elements. In each element $\Omega_e$, the solution is approximated by
\begin{equation}
	\widehat{\uu}_e^{\delta} (\vxi) = \sum_{\bp} \uebp \ell_{\bp} (\vxi)
\end{equation}
where the $\ell_{\bp}$ are tensor-product Lagrange polynomials of degree $N$ which we introduced in~\eqref{eq:lagrange.polynomials} 
to define the map to the reference element~\eqref{eq:ref.map}. The
hat will be used to denote functions written in terms of the reference
coordinates and the delta denotes functions which are discontiuous. Note that
the coefficients $\uebp$ are the unknown values of the function at  solution points
which are taken to be GLL points in this work.
\subsection{Solving for {\bq}}
\begin{figure}
\includegraphics[width=0.95\textwidth]{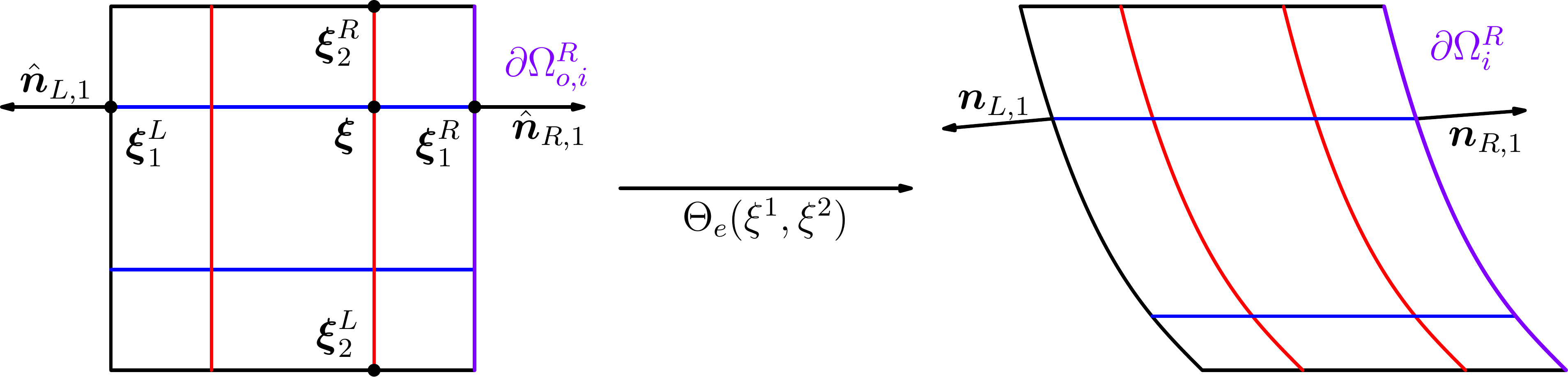}
\caption{\label{fig:ref.map}Illustration of reference map, solution point
	projections, reference and physical normals}
\end{figure}

The system~\eqref{eq:transformed.fo} is solved at each time step for evolving
the numerical solution from time $t^n$ to $t^{n + 1}$, where the first step is solving
the equations for $\bq$. The BR1 scheme was initially introduced for
Discontinuous Galerkin (DG) method in~\cite{Bassi1997} and is used as a
Flux Reconstruction (FR) scheme by exploiting the equivalence between FR and
DG~{\cite{Huynh2007}}, see for example~{\cite{witherden2014}}. 
Here, we show the first step in the FR framework.
Recall that we
defined the multi-index $\bp = (p_i)_{i = 1}^d$ where $p_i \in
\{ 0, 1 \ldots, N \}$. Let $i \in \{ 1, \ldots, d \}$ denote a coordinate
direction and $s \in \{ L, R \}$ so that $(s, i)$ corresponds to the face
$\Ois$ in direction $i$ on side $s$ which has the reference outward normal
$\bnr_{s, i}$, see Figure~\ref{fig:ref.map}. Thus, $\Oip$ denotes the face
where reference outward normal is $\bnr_{R, i} = \be_i$ and $\Oim$ has outward
unit normal $\bnr_{L, i} = - \bnr_{R, i}$. From now, we will suppress the $R$
and denote $\bnr_i = \bnr_{R, i}$.

The FR scheme is a collocation at each of the solution points $\{
\vxi_{\bp} = (\xi_{p_i})_{i = 1}^d, p_i = 0, \ldots, N \}$. We will thus
explain the scheme for a fixed $\vxi = \vxi_{\bp}$ and denote $\vxi_i^s$ as
the projection of {\vxi} to the face $s = L, R$ in the $i^{\tmop{th}}$
direction (see Figure~\ref{fig:ref.map}), i.e.,
\begin{equation}
\left( \vxi_i^s \right)_j = \left\{\begin{array}{lll}
\xi_j, & \qquad & j \neq i\\
- 1, &  & j = i, s = L\\
+ 1, &  & j = i, s = R
\end{array}\right. \label{eq:xis.notation}
\end{equation}
A correction of $\uu_e^{\delta}$ is performed to obtain $\uu_e$ as
\begin{equation}
\uu_e( \vxi ) = \uu_e^{\delta} ( \vxi ) + (
\uu_e^{\ast} - \uu_e^{\delta} ) ( \vxi_i^R ) g_R
(\xi_{p_i}) + ( \uu_e^{\ast} - \uu_e^{\delta} ) ( \vxi_i^L
) g_L (\xi_{p_i}) \label{eq:u.corrected}
\end{equation}
where $\uu^\delta_e(\vxi_i^s)$ denotes the trace value of the normal flux in element $e$ and $\uu_e^{\ast}(\vxi_i^s)$ denotes 
an approximation of the solution at the interface $(s, i)$ which is chosen to be
\begin{equation}
\uu_e^{\ast} ( \vxi_i^s ) = \frac{1}{2}  ( \uu^+_{s, i} + \uu^-_{s, i} ) \label{eq:central.u}
\end{equation}
Thus, $\bq$ can be obtained from~\eqref{eq:transformed.fo} as
\begin{equation}
\bq = \frac{1}{J_{e, p}}  \M \nabla_{\vxi}  \uu_e (\vxi) \label{eq:q.defn}
\end{equation}
\subsection{Time averaging}
The LWFR scheme is obtained by performing the Lax-Wendroff procedure for
Cartesian domains on the transformed equation~\eqref{eq:transformed.fo}. For
$\uu^n$ denoting the solution at time level $t = t_n$ and $\bq^n$ denoting the
gradient computed from~\eqref{eq:q.defn}, the solution at the next time level
can be written as
\[ \uu^{n + 1} = \uu^n + \sum_{k \subindex 1}^N \frac{\Delta t^k}{k!}
\partial_t^{(k)}  \uu^n + O (\Delta t^{N + 2}) \]
where $N$ is the solution polynomial degree. Then, use $\uu_t = - \frac{1}{J}
\nabla_{\vxi} \cdot ( \tfa - \tfv )$ from~\eqref{eq:transformed.fo}
to swap a temporal derivative with a spatial derivative and retaining terms
upto $O (\Delta t^{N + 1})$ to get
\[ 
\uu^{n + 1} = \uu^n - \frac{1}{J}  \sum_{k = 1}^N \frac{\Delta
	t^k}{k!} \partial_t^{(k - 1)}  ( \nabla_{\vxi} \cdot \tfa ) +
\frac{1}{J}  \sum_{k = 1}^N \frac{\Delta t^k}{k!} \partial_t^{(k -
	1)}  ( \nabla_{\vxi} \cdot \tfv ) 
	\]
Shifting indices and writing in a conservative form gives
\begin{equation}
	\uu^{n + 1} = \uu^n - \frac{\Delta t}{J} \nabla_{\vxi} \cdot \tFa +
	\frac{\Delta t}{J} \nabla_{\vxi} \cdot \tFv \label{eq:lw.update}
\end{equation}
where $\tFa, \tFv$ are time averages of the contravariant advective and
viscous fluxes $\tfa, \tfv$
\begin{equation}
\tFs = \sum_{k = 0}^N \frac{\Delta t^k}{(k + 1) !} \partial_t^k  \tfs
\approx \frac{1}{\Delta t}  \int_{t^n}^{t^{n + 1}} \tfs  \ud t, \qquad \alpha = a,v
\label{eq:time.averaged.flux}
\end{equation}
We will find local order $N + 1$ approximations $\tFad_e, \tFvd_e$ to $\tFa,
\tFv$~(Section~\ref{sec:alwp}) which will be discontinuous across element
interfaces. Then, in order to couple the neighbouring elements, continuity of the
normal fluxes at the interfaces will be enforced by constructing the
\tmtextit{continuous flux approximation} using the FR correction functions
$g_L, g_R$~{\cite{Huynh2007}} and the numerical fluxes. Thus, once the local approximations $\tFad_e, \tFvd_e$ are constructed, we construct
the advective and
viscous numerical fluxes for element $e$ in coordinate direction $i \in \{1,\ldots d\}$ at the 
side $s \in \{L,R\}$ (following the notation from~\eqref{eq:xis.notation}) by using Rusanov's~\cite{Rusanov1962} and the central flux respectively
\begin{align}
( \tFa_e \cdot \bnr_{s, i} )^{\ast} & = \frac{1}{2} 
( ( \tFad \cdot \bnr_{s, i} )^+ + ( \tFad \cdot
\bnr_{s, i} )^- ) - \frac{\lambda_{s, i}}{2}  ( \uU^+_{s,
i} - \uU^-_{s, i} ) \label{eq:rusanov.flux.lw} \\
( \tFv_e \cdot \bnr_{s, i} )^{\ast} & = \frac{1}{2} 
( ( \tFvd \cdot \bnr_{s, i} )^+ + ( \tFvd \cdot
\bnr_{s, i} )^- )  \nonumber
\end{align}
The $( \tF^{\delta} \cdot \bnr_{s, i} )^{\pm}$ and $\uU_{s,
	i}^{\pm}$ denote the trace values of the normal flux and time average solution
from outer, inner directions respectively; the inner direction corresponds to
the element $e$ while the outer direction corresponds to its neighbour across
the interface $(s, i)$. The $\lambda_{s, i}$ is a local advection speed estimate at
the interface $(s, i)$. The $g_L, g_R \in \polyP_{N + 1}$ are the FR correction functions~\cite{Huynh2007}, $\uU$ is the approximation of time average solution. The computation of dissipative part of~\eqref{eq:rusanov.flux.lw} using the time averaged solution $\uU$ instead of the solution $\uu$ at time $t_n$ was introduced in~{\cite{babbar2022}} and termed D2 dissipation. It is a natural choice
in approximating the time averaged numerical flux and doesn't add any
significant computational cost because the temporal derivatives of $\uu$ are
already computed during the local approximation of $\tF$. The
choice of D2 dissipation reduces to an upwind scheme in case of constant
advection equation and leads to enhanced Fourier CFL stability
limit~{\cite{babbar2022}}.

The advective and viscous fluxes have been treated separately so far to keep a
simple implementation of the different boundary conditions of the two.
However, for~\eqref{eq:lw.update}, we can combine them to define $\tF =
\tFa - \tFv$, so that the interface numerical fluxes can also be summed as $( \tF_e \cdot \widehat{\bn}_{s, i} )^{\ast} = (
\tFa_e \cdot \widehat{\bn}_{s, i} )^{\ast} - ( \tFv_e \cdot
\widehat{\bn}_{s, i} )^{\ast}$ and thus the continuous time averaged numerical flux is constructed as
\begin{equation}
\begin{split}
( \tF_e ( \vxi ) )^i = ( \tF_e^{\delta} (
\vxi ) )^i & + ( ( \tF_e \cdot \bnr_i )^{\ast} -
\tF_e^{\delta} \cdot \bnr_i ) ( \vxi_i^R ) g_R (\xi_{p_i}) \\
& + ( ( \tF_e \cdot \bnr_i )^{\ast} - \tF_e^{\delta} \cdot
\bnr_i ) ( \vxi_i^L ) g_L (\xi_{p_i})
\end{split}
\label{eq:cts.time.avg.flux}
\end{equation}
Substituting $\tF_e$ in~\eqref{eq:lw.update}, the explicit LWFR update is
\begin{equation}
\begin{split}
\uebp^{n + 1} - \uebp^n + \frac{\Delta t}{J_{e, \bp}}
\nabla_{\vxi} \cdot \tF^{\delta}_e ( \vxi_{\bp} )
+ \frac{\Delta t}{J_{e, \bp}} & \sum_{i = 1}^d ( (
\tF_e \cdot \bnr_i )^{\ast} - \tF^{\delta}_e \cdot \bnr_i )
( \vxi_i^R ) g_R' (\xi_{p_i}) \\ 
+ \frac{\Delta t}{J_{e, \bp}}
& \sum_{i = 1}^d ( ( \tF_e \cdot \bnr_i )^{\ast} -
\tF^{\delta}_e \cdot \bnr_i ) ( \vxi_i^L ) g_L' (\xi_{p_i})
= \bzero
\end{split}\label{eq:lwfr.update.curvilinear}
\end{equation}
\subsubsection{Approximate Lax-Wendroff procedure}\label{sec:alwp}
We now illustrate how to locally approximate $\tFad,
\tFvd$~\eqref{eq:time.averaged.flux} for degree $N = 1$ using the
approximate Lax-Wendroff procedure~{\cite{Zorio2017}}. For $N =
1$,~\eqref{eq:time.averaged.flux} requires $\partial_t  \tfa, \partial_t 
\tfv$ which are
\begin{equation}
\begin{split}
\partial_t  \tfa & \approx \frac{\tfa \left( \uu + \Delta t \uu_t
\right) - \tfa \left( \uu - \Delta t \uu_t \right)}{2 \Delta
t} \backassign \partial_t  \tfad\\
\partial_t  \tfv & \approx \frac{\tfv \left( \uu + \Delta t \uu_t,
\nabla \uu + \Delta t \nabla \uu_t \right) - \tfv \left( \uu -
\Delta t \uu_t, \nabla \uu - \Delta t \nabla \uu_t \right)}{2
\Delta t} \backassign \partial_t  \tfvd
\end{split} \label{eq:ft.defn}
\end{equation}
and $\uu_t, \nabla \uu_t$ are approximated using~\eqref{eq:transformed.fo}
\begin{equation}
\uu_t = - \frac{1}{J} \nabla_{\vxi} \cdot ( \tfad - \tfvd
), \qquad ( \nabla \uu )_t = \frac{1}{J}  \M \nabla_{\vxi}  \uu_t
\label{eq:ut.defn}
\end{equation}
where $\tfad_e, \tfvd_e$ are degree $N$ cell local approximations to the fluxes $\tfa, \tfv$ given
in~\eqref{eq:transformed.fo} constructed by interpolation
\begin{equation}
(\tfad_e)_i  = \sum_{\bp} \tfa (\uu_{e, p})
\ell_{\bp}, \qquad
(\tfvd)_i =\sum_{\bp} \tfv (\uu_{e, p},
\bq_{e, p}) \ell_{\bp}
\label{eq:flux.poly.defn}
\end{equation}
where $\bq$ is obtained in~\eqref{eq:q.defn}. The procedure for other degrees will be similar; Section 4 of~\cite{babbar2022} discusses $N \le 4$ for conservation laws on Cartesian grids.
\subsection{Free stream preservation}
Assume a constant solution $\uu^n = \bc$. The correction terms in~\eqref{eq:u.corrected} will be zero since constant implies continuous at interfaces and thus $\uu_e = \uu_e^\delta = \bc$, proving that $\bq = \bzero$ from~\eqref{eq:q.defn}. Thus, we can now suppress dependence of $\bq$ on $\fv$ and, in particular, write $\tfv = \tfv(\uu_e)$. Thus, the equation we are solving for evolution from $t^n$ to $t^{n+1}$ is now
\[
\uu_t + \frac{1}{J} \nabla_{\vxi} \cdot \tf(\uu) = \bzero, \qquad \tf = \tf^a - \tf^v
\]
implying that the arguments for free stream preservation for LWFR on hyperbolic conservation laws used in~\cite{babbar2024curved} apply to parabolic equations. Thus, as proven in~\cite{babbar2024curved}, the free stream will be preserved as long the metric identity~\eqref{eq:metric.identity.contravariant} is satisfied for interpolated metric terms. See~\cite{Kopriva2006} for methods of approximating the metric terms to discretely satisfy the metric identities.
\section{Numerical results} \label{sec:results}
The numerical experiments were made with the compressible Navier Stokes equations whose description is not provided here to save space, but can be found in~\cite{Gassner_Winters_Hindenlang_Kopriva_2018} or in~\cite{Ray2017} from which some of the numerical tests in this work have been taken. The error based time stepping from~\cite{babbar2024curved} is used in all experiments other than the convergence tests. The extension of error based time stepping of~\cite{babbar2024curved} to parabolic equations is made by using the local time averaged flux to be the total of time averaged advective and viscous fluxes. The subsequent step of truncating the time averaged flux to obtain an embedded lower order method remains the same. Absolute and relative tolerance of $10^{-8}$ is used in all experiments.

The results have been generated using \href{https://github.com/Arpit-Babbar/TrixiLW.jl}{\tt TrixiLW.jl} written using {\tt Trixi.jl}~\cite{Ranocha2022,schlottkelakemper2021purely} as a library. {\tt Trixi.jl} is a high order PDE solver package in {\tt Julia}~\cite{Bezanson2017} and uses the Runge-Kutta Discontinuous Galerkin method. {\tt TrixiLW.jl} has the Lax-Wendroff discretization and uses {\tt Julia}'s multiple dispatch to borrow features like curved meshes support, postprocessing and even the first step of BR1 from {\tt Trixi.jl}. A case for generality of {\tt Trixi.jl} is made since {\tt TrixiLW.jl} is not a fork of {\tt Trixi.jl} but only uses it through {\tt Julia}'s package manager without modifying its internal code. The setup files for the numerical experiments in this work are available at~\cite{nsrepo}.
\subsection{Convergence test}
\begin{figure}
\centering
\begin{tabular}{cc}
\includegraphics[width=0.48\textwidth]{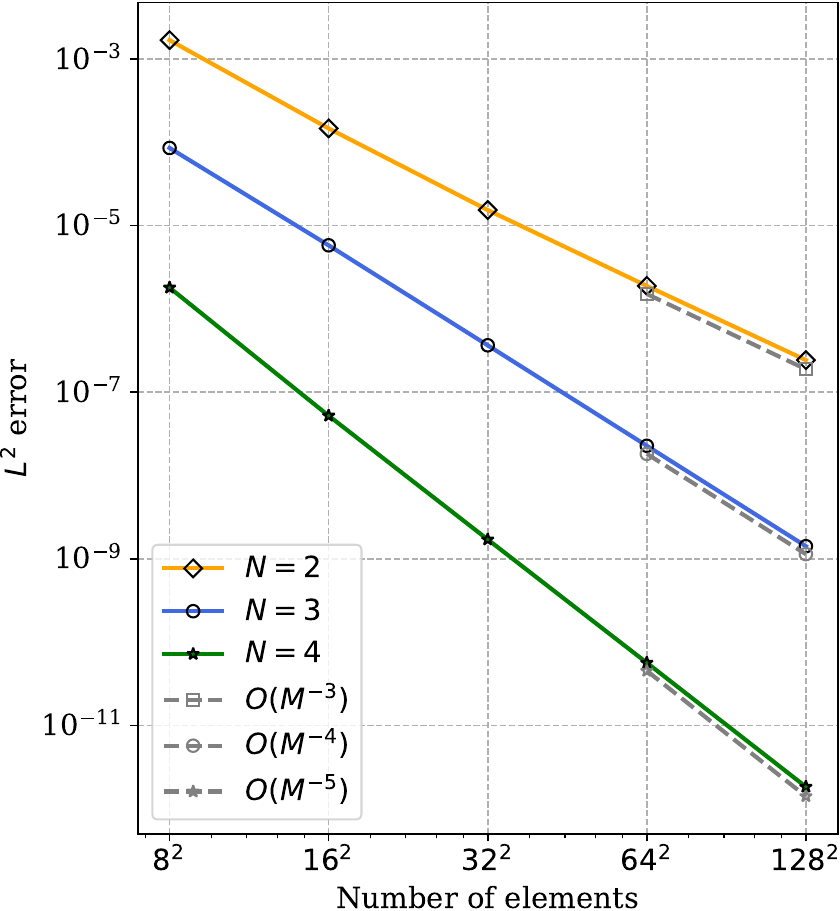} &
\includegraphics[width=0.48\textwidth]{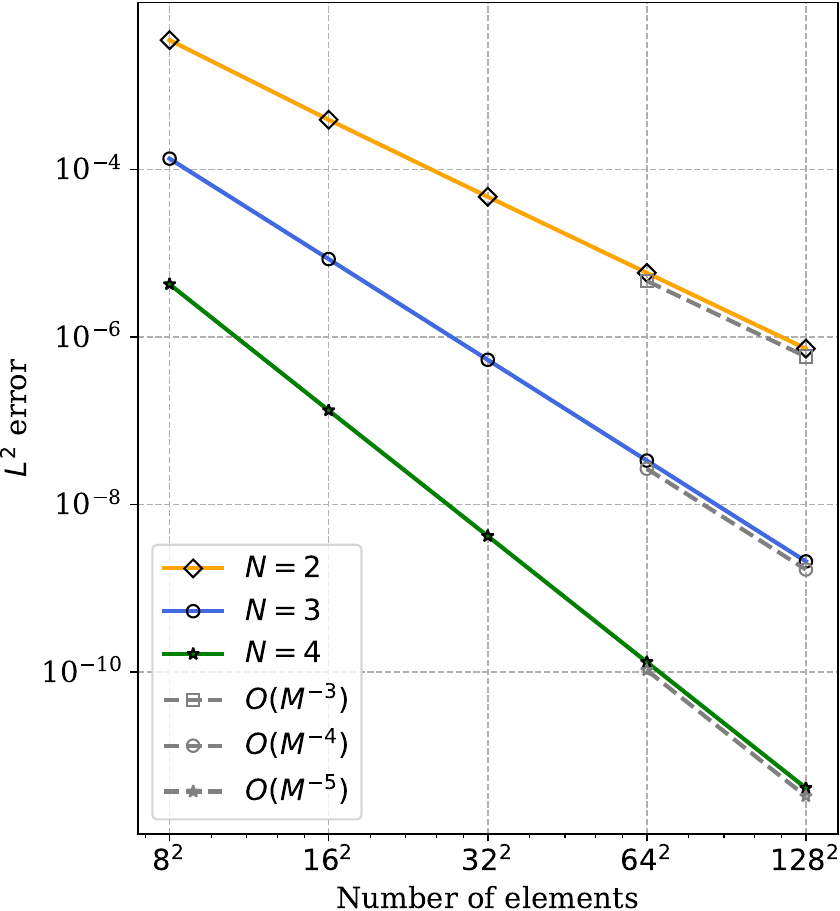} \\
(a) & (b)
\end{tabular}
\caption{Convergence analysis for scalar advection-diffusion equation with $\ba = (1.5,1)$ and coefficient (a) $\nu = 5 \times 10^{-2}$ (b) $\nu = 10^{-12}$\label{fig:convergence.analysis}}
\end{figure}

\begin{figure}
\centering
\begin{tabular}{cc}
\includegraphics[width=0.48\textwidth]{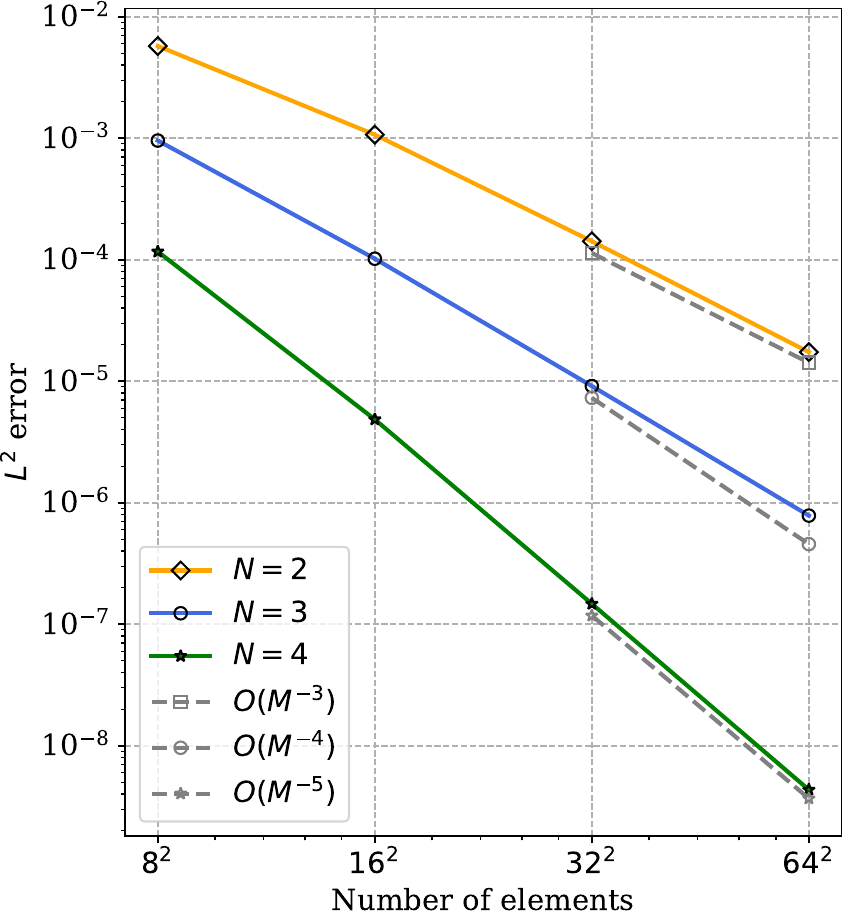} &
\includegraphics[width=0.48\textwidth]{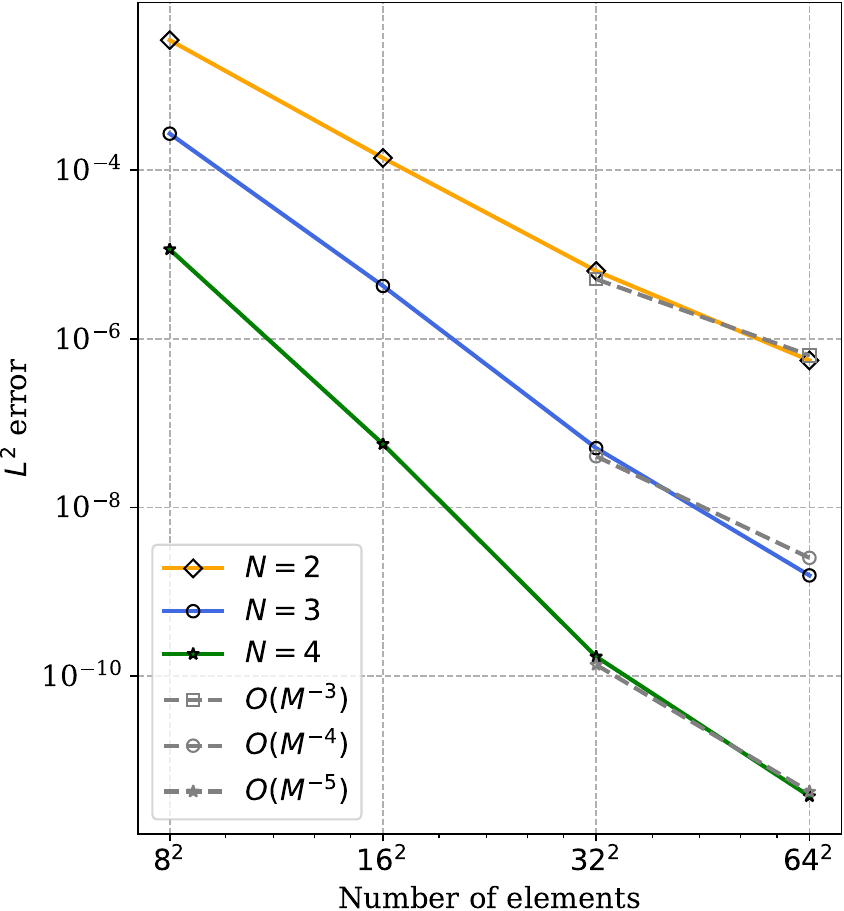} \\
(a) & (b)
\end{tabular}
\caption{Convergence analysis with non-periodic boundary conditions. (a) Eriksson-Johnson test (Section 4 of~\cite{Ellis2016}) and (b) Navier-Stokes equations with manufactured solution}\label{fig:convergence.analysis.nse}
\end{figure}

Consider the scalar advection diffusion equation $u_t + \ba \cdot \nabla u = \nu \Delta u$ where $\ba = (1.5, 1)$. For the initial condition $u_0 (x, y) = 1 + 0.5 \sin(\pi (x + y))$ on $[-1,1]^2$ with periodic boundary conditions, the exact solution is given by $u(x, y) = 1 + 0.5 e^{- 2 \nu \pi^2  t}\sin(\pi(x - a_1 t + y - a_2 t))$. 
A convergence analysis with $\nu$ chosen to be in diffusion and advection dominated regimes is performed and shown in Figure~\ref{fig:convergence.analysis}, and optimal convergence rates are seen for all solution polynomial degrees.
For non-periodic boundaries, we use the Eriksson-Johnson test which is also a scalar advection diffusion whose exact description can be found in Section 4 of~\cite{Ellis2016}. The convergence results are shown in Figure~\ref{fig:convergence.analysis.nse}a where degree $N=2,4$ shown optimal rates while degree $N=3$ nears 3.54 order accuracy. A convergence analysis is also made for the Navier-Stokes equations on the domain $[-1,1]^2$ with a manufactured solution taken from one of the examples in {\tt Trixi.jl}~\cite{Ranocha2022} given by
\begin{equation*}
\begin{split}
\rho &= c + A \sin(\pi x) \cos(\pi y) \cos(\pi t) \\
v_1 &= v_2 = \sin(\pi x) \log(y + 2) (1 - e^{-A(y-2)}) \cos(\pi t) \\
p &= \rho^2
\end{split}
\end{equation*}
The vertical boundaries are periodic and horizontal are no slip, adiabatic walls. The error convergence analysis for density profile is shown in Figure~\ref{fig:convergence.analysis.nse}b, where optimal convergence rates are seen for all polynomial degrees.

\subsection{Lid driven cavity}
\begin{figure}
\centering
\begin{tabular}{cc}
\includegraphics[width=0.45\textwidth]{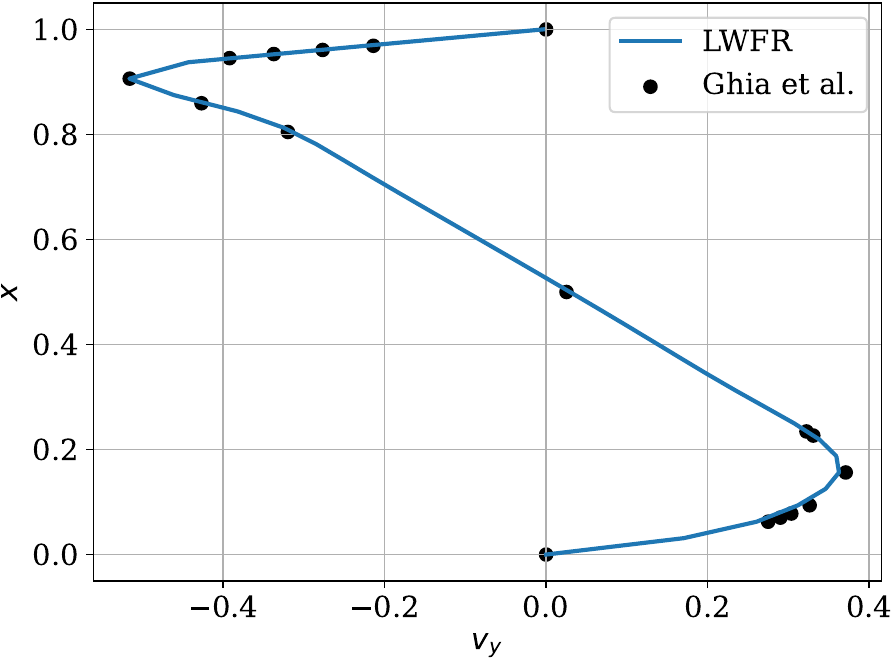}
&
\includegraphics[width=0.45\textwidth]{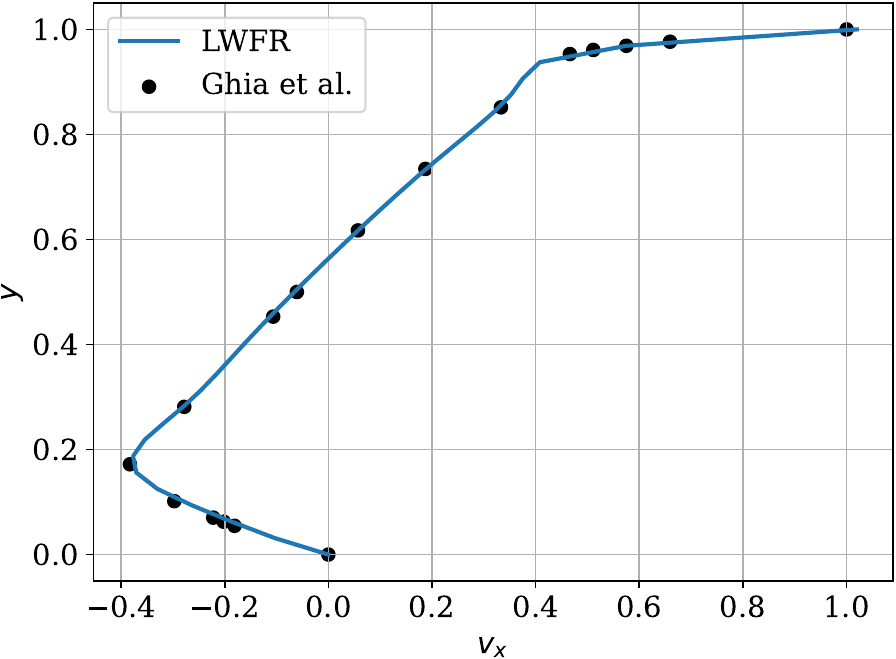}\\
(a) & (b)
\end{tabular}
\caption{Velocity profiles of lid driven cavity test. (a) $v_y$ cut at $y=0.5$ (b) $v_x$ cut at $x=0.5$. \label{fig:cavity}}
\end{figure}

This is a steady state test case for the Navier Stokes equations in the square domain $[0, 1]^2$. We take the setup from~{\cite{Ghia1982}} where {\Pra}=0.7, $\mu = 0.001$ and the initial
condition is the solution at rest
\[
(\rho, u, v, p) = \left(1, 0, 0, \frac{1}{M_{\infty}^2 \gamma} \right)
\]
where $M_{\infty} = 0.1$. The top boundary is moving with a velocity of $(1, 0)$ and the other three have no-slip boundary conditions. All boundaries are isothermal. The problem has a Reynolds number of 1000 corresponding to the top moving wall. The laminar solution of the problem is steady, with Mach number 0.1 corresponding to the moving lid. We compare the solution with the numerical data of Ghia et al.~{\cite{Ghia1982}} by plotting the horizontal velocity profile along the vertical line through the midpoint of the domain, and vertical velocity profile along the horizontal line through the midpoint of the domain. The comparison is shown in Figure~\ref{fig:cavity}, where a good agreement with~{\cite{Ghia1982}} is seen.
\subsection{Transonic flow past NACA-0012 airfoil}\label{sec:swanson}
\begin{figure}
\centering
\begin{tabular}{cc}
\includegraphics[width=0.45\textwidth]{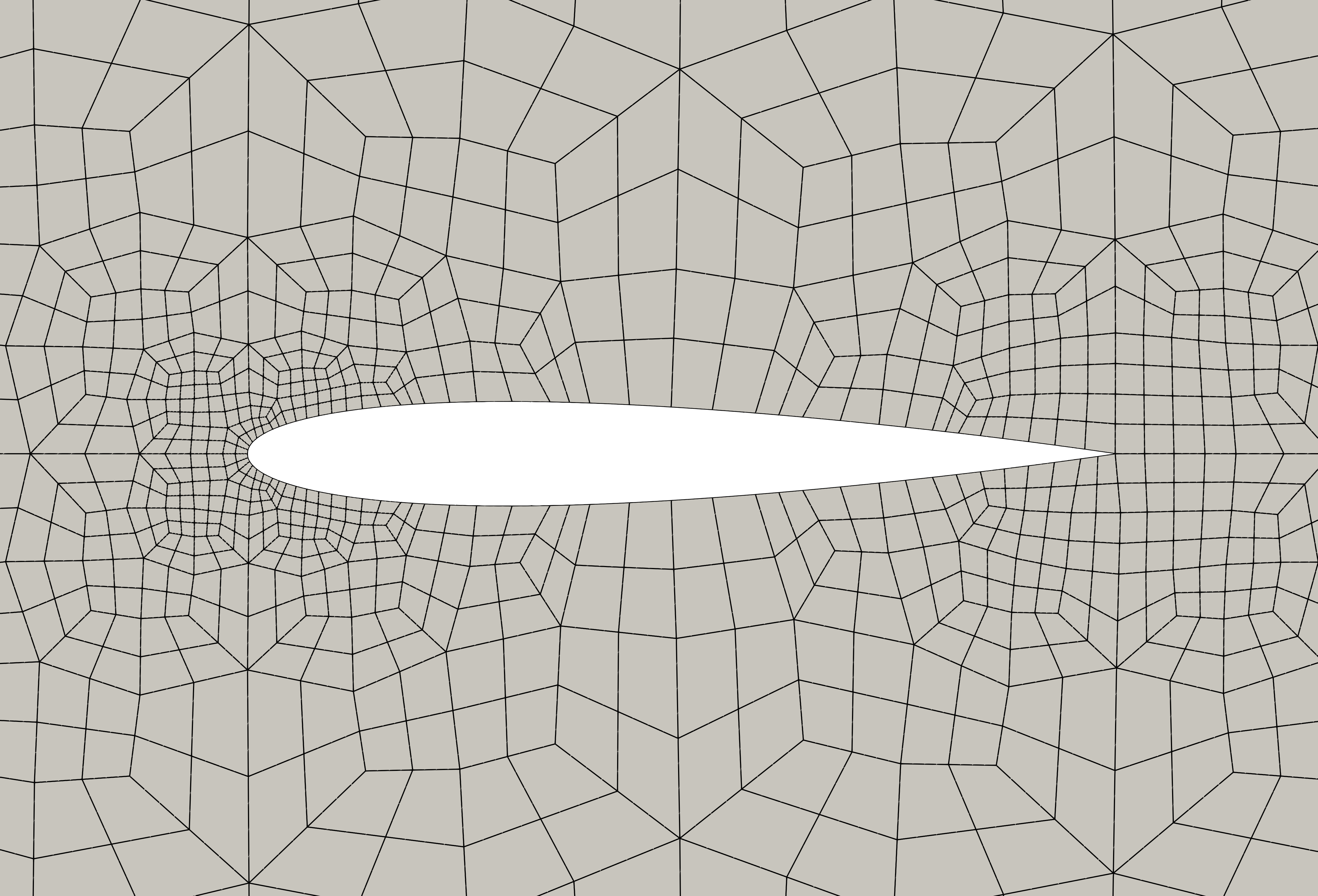} &
\includegraphics[width=0.45\textwidth]{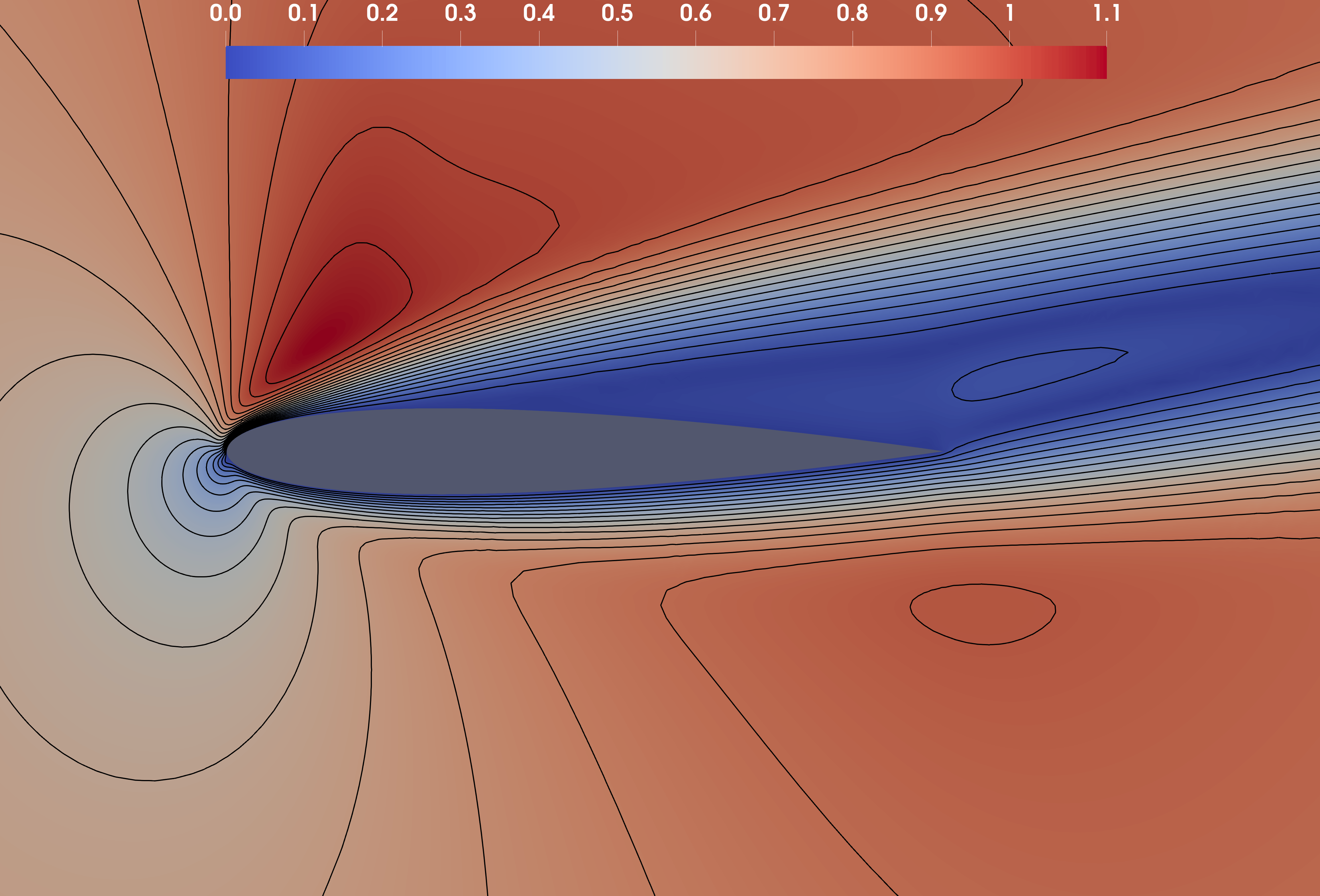}\\
(a) & (b)
\end{tabular}
\caption{Transonic flow over a NACA-0012 airfoil with $M_\infty = 0.8$ solved on a mesh with 728 elements using solution polynomial degree $N=4$. (a) Mesh (b) Mach number contour.}\label{fig:swanson.mesh.mach}
\end{figure}

\begin{figure}
\centering
\begin{tabular}{cc}
\includegraphics[width=0.45\textwidth]{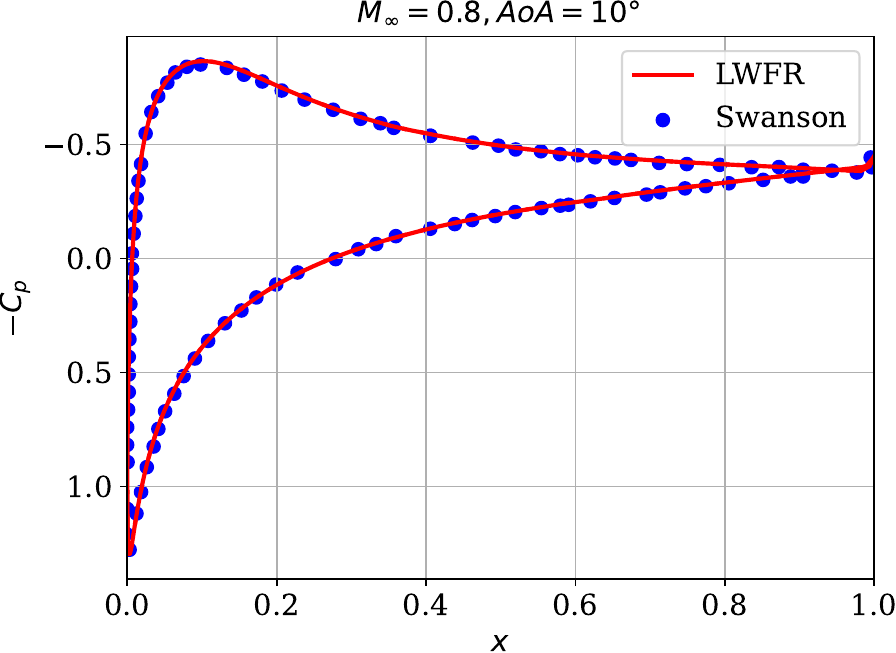} &
\includegraphics[width=0.45\textwidth]{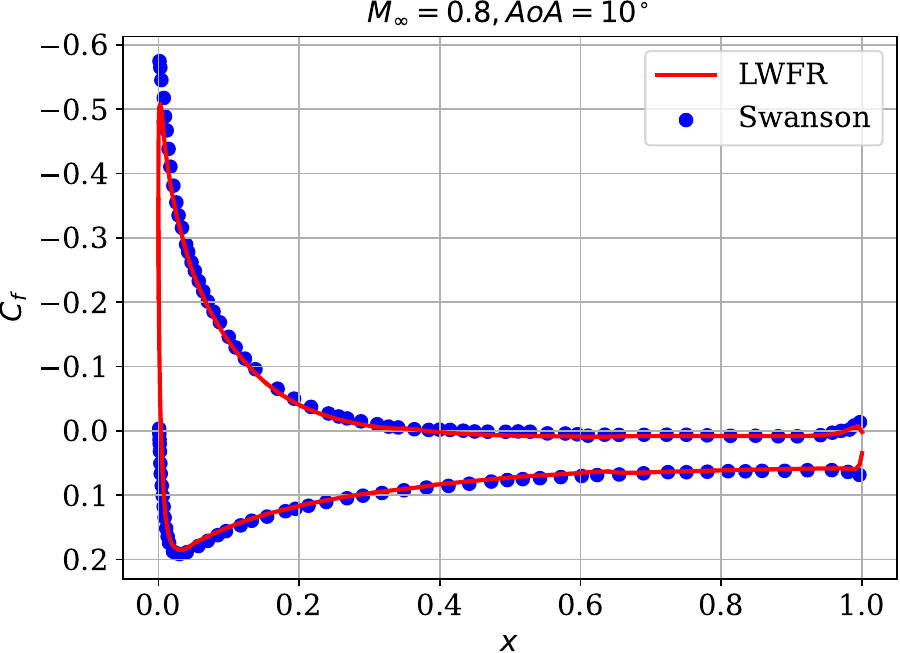}\\
(a) & (b)
\end{tabular}
\caption{Transonic flow over an airfoil, quantities of interest on surface. (a) Coefficient of pressure and (b) Coefficient of friction} \label{fig:cp.cf.swanson}
\end{figure}

This test case involves a steady flow past a symmetric NACA-0012 airfoil. We choose $\rho_\infty = 1, p_\infty = 2.85$ and $Pr = 0.72$, and simulate a flow corresponding to a Reynolds number 500, free-stream Mach number of 0.8 and $10^\circ$ angle of attack. The free-stream velocity is set at $\uu_\infty = (1.574, 0.277)$. The simulation is performed with $728$ elements and polynomial degree $N=4$. The mesh and Mach number contour plot are shown in Figure~\ref{fig:swanson.mesh.mach}. In Figure~\ref{fig:cp.cf.swanson}, coefficient of pressure $C_p$ and coefficient of friction $C_f$ over the airfoil surface are compared with~\cite{Swanson_Langer_2016}, showing good agreement for $C_p$ and same for $C_f$ everywhere other than the leading edge where there are some errors. The coefficients of pressure induced drag and lift ($c_{dp}, c_{lp}$), and drag and lift due to viscous forces ($c_{df}, c_{lf}$) are shown in Table~\ref{tab:swanson.forces} and a good agreement with~\cite{Swanson_Langer_2016} is seen.
\begin{table}
\centering
\begin{tabular}{|c|c|c|c|c|c|}
\hline
& $c_{dp}$ & $c_{df}$ & $c_{lp}$ & $c_{lf}$ & $c_{lp}$+$c_{lf}$\\
\hline
LWFR & 0.1467 & 0.1242 & 0.4416 & -0.0043 & 0.4373\\
\hline
Reference & 0.1475 & 0.1275 & -- & -- & 0.4363\\
\hline
\end{tabular}
\caption{Transonic flow over an airfoil compared with data from~\cite{Swanson_Langer_2016}}\label{tab:swanson.forces}
\end{table}
\subsection{Flow past a cylinder}\label{sec:von.karman}
This test involves a laminar, unsteady flow past a cylinder
inside a channel~\cite{Schafer1996}. The cylinder is slightly offset from the
center of the channel to destabilize the otherwise steady symmetric flow. On
the left, the inflow boundary condition is imposed with $p = \theta =
160.7143$, $v_1 = 4 v_m y (H - y) / H^2$ and $v_2 = 0$, where $v_m = 1.5
\text{m/s}$ is the maximum velocity and the Mach number corresponding to $v_m$
is $0.1$. Isothermal no-slip boundary conditions are imposed on the cylinder
surface, and the top and bottom boundaries. The Reynolds number of the flow
corresponding to the mean velocity is 100, making the viscosity coefficient to be $\nu = 10^{-3}$.

\begin{figure}
\centering
\includegraphics[width=0.85\textwidth]{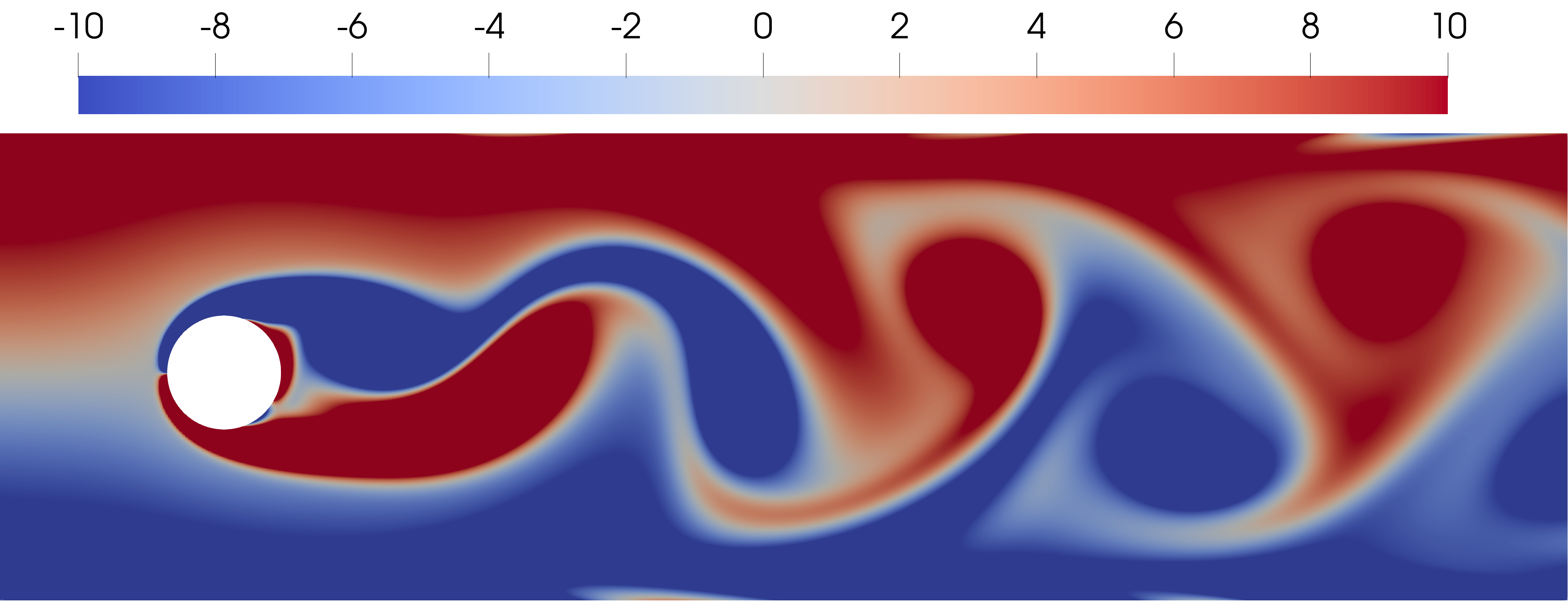}
\caption{Vorticity profile of flow over cylinder}\label{fig:karman.vorticity}
\end{figure}

\begin{figure}
\centering
\begin{tabular}{cc}
\includegraphics[width=0.45\textwidth]{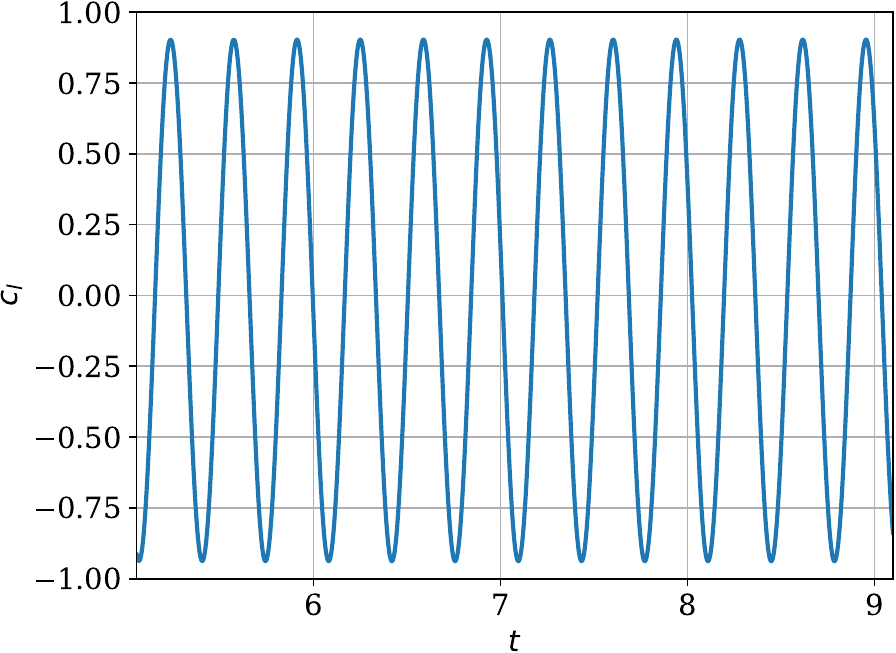} &
\includegraphics[width=0.45\textwidth]{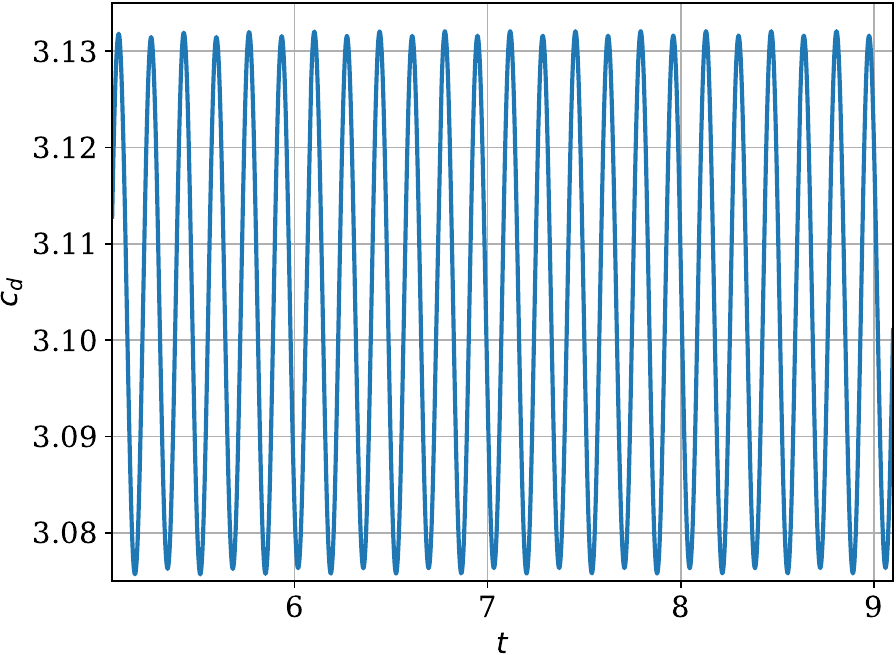}\\
(a) & (b)
\end{tabular}
\caption{Flow over cylinder. (a) Coefficient of lift $c_l$ (b) Coefficient of drag $c_d$}\label{fig:karman.forces}
\end{figure}

The simulation is performed on a nearly uniform mesh with 5692 elements with polynomial degree $N=4$ so that $\Delta x \approx 0.01$. After some time, a Von Karman vortex street appears with a periodic shedding
of eddies from alternate sides of the cylinder. This is typical for slow flows
past a slender body. The vorticity profile is shown in
Figure~\ref{fig:karman.vorticity}, which clearly depicts the periodic vortex
shedding. The periodic behavior can also be observed in
Figure~\ref{fig:karman.forces} where the evolution of the coefficient of total
lift $c_l = c_{l \nocomma p} + c_{l \nocomma f}$ and the coefficient of total
drag $c_d = c_{d \nocomma p} + c_{d \nocomma f}$ on the surface of cylinder at
shown. The time period of the $c_l$ profile is $\tau \approx 0.33759$ so that the
Strouhal number is $\tmop{St} =\mathcal{F}\mathcal{D}/ \overline{u}
=\mathcal{D}/ (\tau \overline{u}) = 0.29621$ where $\mathcal{D}= 0.1,
\overline{u} = 1$ are the diameter of cylinder and mean velocity. This value is in the reference range of~\cite{Schafer1996}. There are some errors in Max $c_l$ and Max $c_l$ as shown in Table~\ref{tab:max.cl.max.cp}, which is to be expected due to the coarse uniform grid which is inefficient for this flow involving multiple boundary layers. The proper resolution of this flow requires appropriately refined grids which will be constructed in future works by adaptive mesh refinement.
\begin{table}
\centering
\begin{tabular}{|c|c|c|c|}
\hline
& Max $c_l$ & Max $c_d$ & St\\
\hline
LWFR & 0.906 & 3.136 & 0.29621\\
\hline
Reference range & $[0.99, 1.01]$ & $[3.22, 3.24]$ & $[0.284, 0.3]$\\
\hline
\end{tabular}
\caption{Comparison of quantities of interest for flow past cylinder}\label{tab:max.cl.max.cp}
\end{table}

\section{Conclusion} \label{sec:conclusion}
The Lax-Wendroff Flux Reconstruction (LWFR) scheme has been extended to
parabolic equations along with its capability of handing curved meshes and
error based time stepping proposed in~{\cite{babbar2024curved}}. The scheme
has been numerically validated by performing convergence and other validation tests 
on the compressible Navier Stokes equations. 
In the sequel, we will develop a shock capturing scheme for adding artificial
dissipation in under-resolved advection dominated problems along with mortar element method for LWFR on parabolic equations.
\bibliographystyle{siam}
\bibliography{references.bib}
\end{document}